\documentclass[12pt]{article}

\usepackage{a4wide}
\usepackage{epsfig}
\usepackage{float}
\usepackage{oldgerm}
\usepackage{tabularx}
\usepackage{xcolor}
\usepackage{amsmath}
\usepackage{placeins}
\usepackage{amssymb}
\usepackage{cite}
\usepackage{multimedia}
\usepackage{amstext,color}
\usepackage{fullpage}

\usepackage{cite}
\usepackage{fixmath}
\usepackage{epstopdf}
\usepackage{MnSymbol}
\usepackage{multirow}
\usepackage{multicol}

\usepackage{graphicx,graphics}
\usepackage{amsthm}
\usepackage[vlined,ruled,linesnumbered]{algorithm2e}

\newtheorem{theorem}{Theorem}[section]

\newtheorem{definition}{Definition}[section]
\newtheorem{lemma}{Lemma}[section]

\def\CC{{\textmd \kern.24em \vrule width.02em height1.4ex depth-.05ex \kern-.26emC}}
\def\TagOnRight
\def\QQ{\rlap {\raise 0.4ex \hbox{$\scriptscriptstyle |$}} {\hskip -0.1em Q}}
\catcode`\@=11
  \def\theequation{\@arabic{\c@section}.\@arabic{\c@equation}} \catcode`\@=12
\begin{document}
\begin{center}
		{{\bf \large {\rm {\bf  Mesh-free mixed finite element approximation  for nonlinear time-fractional biharmonic problem  using weighted $b$-splines}}}}
	\end{center}
 \begin{center}
	{\textmd {\bf Jitesh P. Mandaliya}}\footnote{\it Department of Mathematics,  Institute of Infrastructure, Technology, Research and Management, Ahmedabad, Gujarat, India (jitesh.mandaliya.20pm@iitram.ac.in)},
	{\textmd {{\bf Dileep Kumar}}}\footnote{\it Department of Mathematics, Government Post Graduate College Noida, Uttar Pradesh, India (dilipkmr832@gmail.com)},
 {\textmd {{\bf Sudhakar Chaudhary}}}\footnote{\it Department of Mathematics,  Institute of Infrastructure, Technology, Research and Management, Ahmedabad, Gujarat, India (dr.sudhakarchaudhary@gmail.com)}
\end{center} 
 \begin{abstract} In this article, we propose a fully-discrete scheme for the numerical solution of  a nonlinear time-fractional biharmonic problem. This problem is first converted into an equivalent system by introducing a new variable. Then spatial and temporal discretizations are done by the weighted $b$-spline method and $L2$-$1_\sigma$ approximation, respectively. The weighted $b$-spline method uses weighted $b$-splines on a tensor product grid as basis functions for the finite element space and by construction, it is a mesh-free method. This method combines the computational benefits of $b$-splines and standard mesh-based elements. We derive $\alpha$-robust   \emph{a priori} bound and convergence estimate in the $L^2(\Omega)$ norm for the proposed scheme.  Finally, we carry out few numerical experiments to support our theoretical findings. 
\end{abstract}
	\section{Introduction}
	\par For a bounded convex polygonal domain $\Omega \subset  \mathbb{R}^{2}$, we consider the following time-fractional biharmonic problem:
	\begin{subequations}
		\label{bha1}
		\begin{align}
			\begin{split}
				D_t^\alpha u+\Delta^2 u-\Delta u=& \  f(u)\ \ \mbox{in} \   \Sigma,
			\end{split}
			\label{bha2}\\
			\begin{split}
				u=\Delta u=& \ 0 \quad \ \ \ \mbox{on}\ \partial\Sigma,
			\end{split}
			\label{bha3}\\
			\begin{split}
				u({x},0)=& \ u_{0}({x}) \ \mbox{in} \; \Omega,
			\end{split}
			\label{bha4}
		\end{align}
	\end{subequations}
	 where $\partial \Omega$  is the domain boundary, $\Sigma=\Omega\times (0,T]$, $\partial\Sigma=\partial\Omega\times (0,T]$,  $T>0$ and $f(u)=u-u^3$. In \eqref{bha1}, $D_t^\alpha u$ represents the  $\alpha^{th}$ order Caputo derivative of $u$   which is defined below \cite{[r1]}: \\
	\begin{equation}\label{cc}
		{D^{\alpha}_{t}}u(x,t):= \frac{1}{\Gamma (1-\alpha)}\int_{0}^{t}(t-s)^{-\alpha} \; \frac{\partial u(x,s)}{\partial s} \, ds \quad \text{for} \;\;  \ \  0 < \alpha < 1.
	\end{equation}
In recent years many authors have studied the fourth order partial differential equations (PDEs) (see \cite{[yliu],[amiya1],[nrcs],[kaara],[hung],[fb1],[kaara1],[vfem],[gudi]}  and references therein). 
  Authors in \cite{[amiya1]} presented a mixed finite element method (FEM) with $C^0$-piecewise linear elements for solving a nonlinear parabolic biharmonic problem (integer-order time derivative)  with Navier boundary condition. In \cite{[fb1], [hung]}, a linear time-fractional biharmonic problem is explored and  $L1$ scheme along with mixed FEM is used to find the numerical solution. Authors in \cite{[kaara]} have discussed the existence-uniqueness of problem \eqref{bha1} and used the FEM in space and backward Euler convolution quadrature in time to obtain the numerical solution. They also derived optimal order error estimates for smooth and non-smooth data.  In \cite{[yliu]}, authors considered a nonlinear time-fractional biharmonic problem. To obtain a numerical solution, they utilized a two-grid algorithm that consist of two distinct steps: first, they solved a nonlinear system on a coarse grid using nonlinear iterations then they solved the linearized system on the fine grid by Newton iterations. A scheme based on the orthogonal spline collocation method and $L1$ approximation formula is proposed in \cite{[nbha]} for solving a nonlinear fourth-order reaction–subdiffusion equation. In \cite{[l21]}, authors explored a high-order two-grid compact difference method for a nonlinear time-fractional biharmonic equation using the well-known $L2$-$1_{\sigma}$ scheme. The solution of time-fractional PDEs in general exhibits weak singularity at $t=0$. For the time-fractional biharmonic problem, in \cite{[fb1], [kaara], [hung], [nbha]}, the weak singularity at $t=0$ of the unknown solution has been taken into consideration and graded mesh is used for the temporal discretization.
 \par  Geometry description and mesh generation are generally challenging and time-consuming tasks in FEM. H$\ddot{o}$llig et al. in \cite{[HRW00a]} introduced a mesh-free method called the weighted extended $b$-spline (web-spline) method. The web-spline method involves weight function \cite{[HRW00a],[rvach]}, extension procedure \cite{[HRW00a],[kh]}, and 
 $b$-splines \cite{[de2]}. Several research articles in the literature have explored the applications of the web-spline method to compute the numerical solution for the elliptic \cite{[nrcs], [sv1], [HRW00a]} and parabolic \cite{[JC], [ptv]} PDEs with integer-order derivatives.  In \cite{[nrcs]}, authors studied the following steady state biharmonic problem: 
 \begin{subequations}
		\label{eebha1}
		\begin{align}
			\begin{split}
				  \Delta^2 u& =f \quad \mbox{in} \ \  \Omega,
			\end{split}
			\label{eebha2}\\
			\begin{split}
				u= \Delta u&=  0, \quad \mbox{on}\ \partial \Omega. 
			\end{split}
		\end{align}
	\end{subequations}
 By substituting $v=-\Delta u$, they split problem \eqref{eebha1} into two second order problems and then used the web-spline mesh-free method to get the numerical solution. For problem \eqref{bha1}, this splitting is not feasible because our problem is more general due to the presence of Caputo derivative \cite{[fb1]}. 
\par Motivated by the literature mentioned above, in this work, we have employed the weighted $b$-spline method with the $L2$-$1_{\sigma}$ scheme to solve the time-fractional biharmonic problem \eqref{bha1}. Note that the weighted $b$-spline method is basically web-spline method without extension procedure. The main contributions of the present work are given below:\\
 (i) We propose a mesh-free method based on $L2$-$1_{\sigma}$ approximation and weighted $b$-splines. This method gives smooth and high-order accurate approximations in spatial direction with relatively few parameters (degree of freedom). \\
 (ii) In the weighted $b$-spline method,  essential boundary conditions can be imposed exactly with the help of an appropriate choice of weight function.\\ 
 (iii) The $L1$ scheme with mixed FEM is analyzed in \cite{[fb1]} for the linear time-fractional biharmonic problem. In comparison with $L1$ scheme, the $L2$-$1_{\sigma}$ scheme exhibits superior convergence rate. However, the analysis of problem \eqref{bha1} using  $L2$-$1_{\sigma}$ based mixed FEM, introduces more challenges than its $L1$ counterpart. Here, we derive the error estimates for a $L2$-$1_\sigma$ based fully-discrete scheme. These estimates are $\alpha$-robust (bounds do not blow up as $\alpha$ $\rightarrow$ $1^{-}$).\\
 (iv) To illustrate the advantages of the proposed method, we perform some numerical experiments for weighted $b$-splines of degree $\le 3$. We also compare the results of the proposed method with mixed FEM for linear basis functions.\\
 To the best of our knowledge, this is the first attempt to use a weighted $b$-spline based mesh-free method to solve the problem \eqref{bha1}.
 \par The paper is organized as follows: In Section 2, we construct the weighted $b$-spline basis for spatial discretization. In Section 3, a fully-discrete scheme is introduced for solving the problem \eqref{bha1}. This scheme utilizes the $L2$-$1_\sigma$ approximation on a graded mesh in temporal direction and weighted $b$-splines in spatial direction. An $\alpha$-robust   \emph{a priori} bound and convergence result are derived for the numerical solution within the framework of the proposed scheme. In Section 4, we present some numerical experiments that validate our theoretical convergence results. Finally, Section 5 presents the conclusions drawn from our study.
		\section{Weighted $B$-Spline Basis}
	\noindent In this section, we describe the construction of the weighted $b$-spline basis functions \cite{[HH],[kh]}. This construction comprises of two essential components: the weight function and the $b$-splines. First, we define the $b$-splines.

	 \noindent For a given knot sequence $\xi:\, $\dots$ \le \xi_0\le\xi_1\le\xi_2\le \,\dots$,  the $b$-splines $b_{k,\xi}^{m}$ of degree $m$ are defined by the following recursion formula \cite{[de2]}:
	
	\begin{equation}\label{bb1}
			b_{k,\xi}^{m}(x) = \frac{x-\xi_{k}}{\xi_{k+m}-\xi_{k}}b_{k,\xi}^{m-1}(x)+\bigg(1-\frac{x-\xi_{k+1}}{\xi_{k+m}-\xi_{k+1}}\bigg)b_{k+1,\xi}^{m-1}(x), 
		\end{equation} 
	where $x \,\in\, \mathbb{R}$,   $k\, \in \, \mathbb{Z}$, 	
	\begin{equation*}
			b_{k,\xi}^{0}(x)= \, \left\{\begin{array}{l} 1 \qquad \mbox{for} \quad \xi_{k} \leq x < \xi_{k+1}, \\ 0  \qquad \mbox{otherwise}, 
				\end{array} \right.
		\end{equation*}
	and discarding terms for which the denominator vanishes.
	\noindent For a knot sequence $0, 1,\dots, m+1,$ the recursion \eqref{bb1} reduces to the following form: 
	\begin{equation*}
		b^m(x)=\frac{x}{m}b^{m-1}(x)+\frac{m+1-x}{m}b^{m-1}(x-1).
	\end{equation*}
		\noindent The $b^m$ is also known as a standard uniform $b$-spline of degree $m$. A $d$-variate uniform $b$-spline $b_{k}$ of degree $m$ and  grid width $h$ is a product of scaled translates of univariate cardinal $b$-splines \cite{[HH]}:
	\begin{equation*}
		b_k(x)= \prod_{i=1}^{d} b^m\left({\frac{x_i}{h}}-k_i\right),
	\end{equation*}
	where $k=(k_1, k_2,\dots, k_d) \in \mathbb{Z}^d$,  $x=(x_1, x_2,\dots, x_d)\, \in \mathbb{R}^d$ and $h$ is the grid width of cell $Q_l=lh+[0,1]^dh, \ l\in \mathbb{Z}^d$. 
\\ \\
\begin{multicols}{2}
\begin{figure}[H]
	\centering
	\includegraphics[width=8cm]{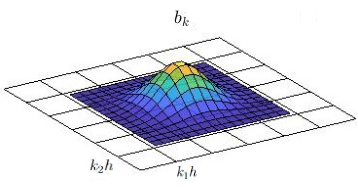}
	\caption{\emph{Uniform $b$-spline}}
	\label{picb}
\end{figure}
\columnbreak
\begin{figure}[H]
	\centering
	\includegraphics[width=5cm]{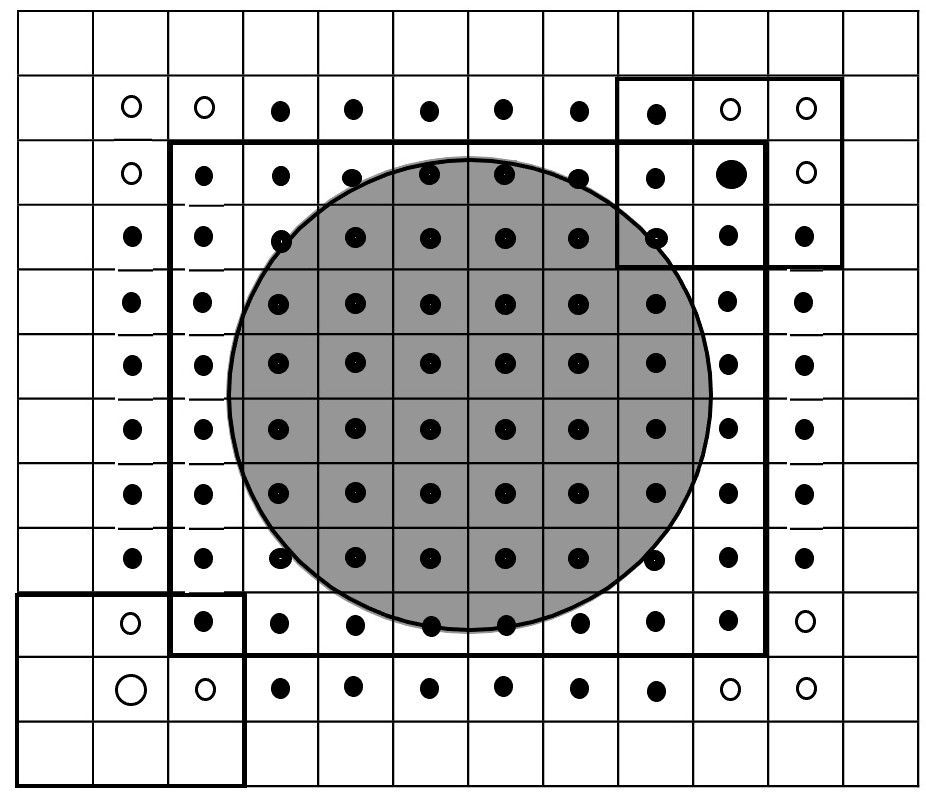}
	\caption{\emph{$B$-spline on bounded domain}}
	\label{picb1}
\end{figure}
\end{multicols}
 As is illustrated in Figure \ref{picb}, in general $b_k$ is a nonnegative, $(m-1)$-times continuously differentiable in each variable with support $\mbox{supp} \, b_k=kh+[0, m+1]^dh$. On each grid cell $Q_l$,  $b$-spline $b_k$ is a polynomial of degree $m$ in each variable. Those $b$-splines that have some support in $\Omega$ are called relevant $b$-splines otherwise, they are irrelevant $b$-splines. In Figure \ref{picb1}, relevant and irrelevant biquadratic $b$-splines are marked with dots and circles at the centers of their supports, respectively. For any open set $\Omega \subset \mathbb{R}^d$, the relevant $b$-splines are linearly independent \cite{[colweb]}. Let $V_h$ be the space spanned by all relevant $b$-splines, then we call  $V_h$ as spline space and write $V_h=span\{b_k: k\in K\},$ where $K$ is the relevant index set. In general, $b$-splines are not utilized as a finite element space due to two reasons. Firstly, $b$-splines do not conform to essential boundary conditions, and secondly, they do not provide a stable basis. To overcome these issues, H$\ddot{o}$llig et al. in \cite{[HRW00a],[kh]} introduced a weight function and extension procedure. 
	
	Weight functions can be determined analytically for specific domains, such as circles, squares, etc. For more details on constructing the weight functions $w$, various methods are available as discussed in \cite{[kh],[HRW00a],[rvach]}.

The stability issue of $b$-spline basis functions has been addressed by H$\ddot{o}$llig et al. \cite{[HRW00a]} by utilizing the extension procedure. Nonetheless, an alternative perspective is presented by authors in \cite{[HH]}, demonstrating that accurate approximation can be achieved in practical computations even without employing the extension procedure. Notably, they emphasized that basic preconditioning techniques suffice for solving ill-conditioned finite element systems with acceptable precision.
	Now, we define weighted $b$-spline space as follows \textup{\cite{[HH]}}:
	\begin{equation*}
		B_h= span \{wb_k: k \in K\},\\
	\end{equation*}
	where $\{b_k\}_{k \in K}$ are  $b$-splines of degree $m \geq 1$.
	 We also consider irrelevant indices for computation. This implies that we can write the approximate solution $u_h$ in the following form: \[u_h=\sum_{k \in K^{'}} \beta_kb_k,\] where $K^{'}$ represents the smallest rectangular array containing $K$, and $\beta_k$ are set to zero for $k\notin K$.\\
	 Hereafter $C>0$  denotes a generic constant that can take different values
at different occurrences and depending on several parameters like domain $\Omega$, degree of $b$-spline $m$, and the weight function $w$, etc. However, it is always independent of spatial and temporal step sizes.	
	\section{Fully-discrete formulation}
 \noindent For the convex domain, we can rewrite problem \eqref{bha1} into following  system of two second order equations: 
	\begin{equation}\label{mixed}
			\begin{split}
				{D_t^\alpha} u-\Delta v+v&= \  f(u)\  \qquad    \mbox{in}  \qquad \Sigma,\\
				-\Delta  {u}&=v \qquad \quad \ \ \ \,\mbox{in} \qquad \Sigma, \\
				{u}({x},0)&= \ {u}_{0}({x}) \qquad  \mbox{in } \ \ \quad \Omega,\\
				u= &v= \ 0 \qquad \quad  \ \mbox{on} \qquad   \partial \Sigma.
			\end{split}
	\end{equation}
	\begin{definition}
		Weak formulation of \eqref{mixed} is to find $u(\cdot,t)\  \mbox{and}\ {v(\cdot,t)}\in  H_{0}^{1}(\Omega)$ such that for each $t\in (0,T]$
			\begin{align*}
				\begin{split}
					\left(D_t^\alpha u,\phi\right)+(\nabla v,\nabla \phi)+( v,\phi)&=(f(u),\phi) \ \quad \forall \phi \in H_{0}^{1}(\Omega),
				\end{split}\\
				\begin{split}
					(\nabla u,\nabla \psi)&=( v, \psi) \qquad    \quad\forall \psi \in H_{0}^{1}(\Omega),
				\end{split}\\
				\begin{split}
					u(x, 0)&= u_{0}(x) \qquad \quad \   \forall x \in \Omega.
				\end{split}
			\end{align*}
	\end{definition}
		For fully-discrete formulation, we divide time interval $[0, T]$ using $t_n = T \left(\frac{n}{N}\right)^r$ for $n = 0, 1, . . . , N$, where the grading parameter $r\geq 1$. For each $n= 0, 1,..., N$, let $u^n := u(t_n)$ and $U_h^n$ denote the approximate value of $u$ at $t_n$. We approximate Caputo fractional derivative by $L2$-$1_\sigma$ approximation \cite{[AAl2],[r16]}: \\ 
	\begin{equation}\label{cfd}
		\begin{split}
			D^{\alpha}_{ t_{n-\sigma}}u \, \approx& \, (D^{\alpha}_{N} u)^{n-\sigma}\\
			:=& \sum_{j=1}^{n}\mathcal{A}_{\sigma,n-j}^{n} (u^{j}-u^{j-1}) \quad \mbox{for} \ n=1, 2,\dots, N.  
		\end{split}
	\end{equation}
 For $\sigma=\frac{\alpha}{2}$ and $t_{n-\sigma}=\sigma t_{n-1}+(1-\sigma)t_n$, the coefficient $ \mathcal{A}_{\sigma,n-j}^{n}$ in \eqref{cfd} is given by $ \mathcal{A}_{\sigma,0}^{1}=a_{1,0}$ if $n=1,$ and for $n\geq 2$  
	$$\mathcal{A}_{\sigma,n-j}^{n} = \begin{cases}  a_{n,0}+\frac{\tau_{n-1}}{\tau_n}b_{n,1} & \mbox{} \  \  \mbox{if}\  j=n, \\ a_{n,n-j}+\frac{\tau_{j-1}}{\tau_j}b_{n,n-j+1}-b_{n,n-j} & \mbox{ } \mbox{ if} \  2 \leq j\leq n-1, \\ a_{n,n-1}-b_{n,n-1} & \mbox{ } \  \mbox{if} \ j=1,\\\end{cases}$$
	where
	\begin{equation*}
		\begin{split}
			a_{n,0}&=\frac{\tau_n^{-1}}{\Gamma(1-\alpha)}\int_{t_{n-1}}^{t_{n-\sigma}}(t_{n-\sigma}-\eta)^{-\alpha}d\eta, \ \\
			a_{n,n-j}&=\frac{\tau_j^{-1}}{\Gamma(1-\alpha)}\int_{t_{j-1}}^{t_{j}}(t_{n-\sigma}-\eta)^{-\alpha}d\eta \quad  \ \mbox{for}   \ \ 1\leq j\leq n-1,\\
			b_{n,n-j}&=\frac{2\tau_j^{-1}}{\Gamma(1-\alpha)(t_{j+1}-t_{j-1})}\int_{t_{j-1}}^{t_{j}}(t_{n-\sigma}-\eta)^{-\alpha}(\eta-t_{j-\frac{1}{2}})d\eta \  \  \mbox{for} \  1\leq j\leq n-1.
		\end{split}
	\end{equation*}
	For our further analysis, we set some notations that are:
	\[U_h^{n,\sigma} := \sigma \, U_h^{n-1} + (1- \sigma) \, U_h^{n}, \quad V_h^{n,\sigma} := \sigma \, V_h^{n-1} + (1- \sigma) \, V_h^{n}, \ \ \quad \ \ \  \mbox{for} \ \ n\geq1.\] \[\tau_{n,\sigma}^*:=(1-\sigma)\frac{\tau_n}{\tau_{n-1}}, \quad \mbox{and} \quad \hat{U}_h^{n-1,\sigma} := (1+\tau_{n,\sigma}^*) \, U_h^{n-1} - \tau_{n,\sigma}^* \, U_h^{n-2} \ \  \mbox{for} \ \ n\geq2.\]
	\noindent Then the fully-discrete scheme is to compute $U_h^n$ and $V_h^n\in B_h$ such that  $U_h^0=u_h^0$ for $n=0$ and 
\begin{subequations}
		\label{fd}
		\begin{align}
			\begin{split}
				(({D_N^\alpha}U_h)^{1-\sigma},x_h)&+(\nabla V_h^{1,\sigma},\nabla x_h)+( V_h^{1,\sigma}, x_h)\\
				&=(f({U}_h^{1,\sigma} ),x_h) \  \forall x_h\in B_h,
			\end{split}
			\label{fnd1}\\
			\begin{split}
				(({D_N^\alpha}U_h)^{n-\sigma},y_h)&+(\nabla V_h^{n,\sigma},\nabla y_h)+( V_h^{n,\sigma}, y_h)\\
				&=(f(\hat{U}_h^{n-1,\sigma}), y_h) \  \forall y_h\in B_h \ \mbox{for} \ \ n\geq 2,
			\end{split}
			\label{fd1}\\
			\begin{split}
				(\nabla U_h^{n},\nabla z_h)&=( V_h^{n}, z_h) \quad \forall z_h\in B_h \ \mbox{for}  \ n\geq 1,
			\end{split}
			\label{fd2}
		\end{align}
	\end{subequations}
 
where $u_h^0$  is some approximation of $u_0.$\\
 Now using mathematical induction and inverse inequality as in \cite{[Dli]}, one can find the bound for the solution $U_h^{n-1}$ in $L^{\infty}$ norm as
 \begin{equation}\label{ll}
     \|U_h^{n-1}\|_{L^{\infty}} \le \|R_h u^{n-1}\|_{L^{\infty}}+\|R_h u^{n-1}-U_h^{n-1}\|_{L^{\infty}}\le C  \ \mbox{for}  \ n\geq 1. 
 \end{equation}
 Use of \eqref{ll} yields the following estimate for nonlinear term $f$ given by
 \begin{equation}\label{lip}
     \|f(u^{n-1})-f(U_h^{n-1})\| \le L_f\|u^{n-1}-U_h^{n-1}\| \ \mbox{for}  \ n\geq 1.
 \end{equation}
 \section{ { Weighted $b$-splines Based Estimates }}
 Here we introduce $L^2(\Omega),$ Ritz projectors, discrete Laplacian operator, and their relationship \cite{[hr12]} that will be required in the derivations of a priori bound and convergence results.\\
 Now, we define $L^2$ projector often denoted by $P_h$ is a map from $L^2(\Omega)$ into $B_h$ such that 
 \begin{equation}\label{DK25}
			(w, v_h)=(P_h w,v_h)\quad \forall\   v_h \in B_h,
	\end{equation}
the Ritz projector denoted by $R_h$ is a map from $H_0^1(\Omega)$ into $ B_h$ such that
\begin{equation*}
		\begin{split}
			(\nabla w,\nabla v_h)=&(\nabla R_h w,\nabla v_h)\quad \forall \ v_h \in B_h,
		\end{split}
	\end{equation*}
 the discrete Laplacian operator  is denote by $\Delta_h$ is a map from $B_h$ into $B_h$ such that   
	\begin{equation*}
		(\nabla w,\nabla v_h)=(-\Delta_h w,v_h) \quad \forall \ v_h \in B_h,
	\end{equation*}
 and these operators are related to each other by the following relation 
 \begin{equation}\label{pr}
    \Delta_hR_hw=P_h\Delta w\quad \forall\  w\in H^2(\Omega).
 \end{equation}
	Here, we also require the discrete coefficients for the further analysis that are denoted by $Q^{(n)}_{n-i}$ for $n=1, 2,\dots, N$ and defined by
	\begin{equation*}
		\begin{split}
			Q^{(n)}_{n-i} :=  \left\{\begin{array}{l} \frac{1}{\mathcal{A}_{\sigma,0}^{i}}\sum_{k=i+1}^{n} \left(\mathcal{A}_{\sigma,k-i-1}^{k} - \mathcal{A}_{\sigma,k-i}^{k} \right) Q^{(n)}_{n-k} \quad \mbox{if} \; \  i=1, 2,\dots, n-1, \\
				\\
				\frac{1}{\mathcal{A}_{\sigma,0}^{n}} \quad \mbox{if} \; \, i=n. \\
			\end{array}\right.
		\end{split}
\end{equation*}
 Furthermore, these discrete coefficients $Q^{(n)}_{n-i}$ satisfy the following estimate for each $n=1, 2,\dots, N$.
 \begin{lemma}\label{l7}
		\cite{[r15]} If
		the constant $\gamma \in (0,1)$ then for each  $n=1, 2,\dots, N$ one has

		\begin{equation*}
				\sum_{j=1}^{n} \, Q^{(n)}_{n-j} \, j^{r (\gamma - \alpha)} \, \le \, \frac{\, 11\Gamma(1 + \gamma - \alpha)}{ \, 4\Gamma(1+ \gamma)} \, T^{\alpha} \, \bigg(\frac{t_n}{T}\bigg)^{\gamma} \, N^{r (\gamma - \alpha)}.   
		\end{equation*}	
	\end{lemma} 
 Now, we state two important lemmas called coercivity property and the discrete fractional Gronwall inequality that help us to derive $\alpha$-robust estimates for the fully-discrete solution $(U_h^n,V_h^n)$.
 \begin{lemma}\label{l111}
		\cite{[AAl2]} If $\{\omega^k\}_{k=1}^{N}$ is an arbitrary sequence of functions in $L^2(\Omega).$ Then it satisfies
		\begin{equation*}
			\left( ({D}^{\alpha}_{N} \omega)^{n-\sigma}, \, \omega^{n,\sigma} \right) \, \ge \, \frac{1}{2} \, ({D}^{\alpha}_{N} \|\omega\|^2)^{n-\sigma}    \ \text{for} \ n=1, 2,\dots, N.
		\end{equation*}
	\end{lemma}
	\begin{lemma}\label{l5}
		\cite{[r16]} For the constants  $0 \le \gamma \le 1,$ $\lambda_l\ge 0,$ and $\Lambda$ with $ \sum_{l=0}^{n} \lambda_l \le \Lambda$ for $n\geq 1$. If the sequences $\{\zeta ^n\ge 0, \: \Theta ^n\ge 0\ \text{for} \ n \ge 1\}$ are bounded with the sequence $\{\omega^n\ge 0 \ \text{for} \ n \ge 0\}$ satisfies
		\begin{equation*}
			\begin{split}
				(D_N^{\alpha} \omega^2)^{n- \gamma} \, \le \, \sum_{k=1}^{n} \, \lambda_{n-k} \, (\omega^{k, \theta_1})^2 + \, \zeta^n \, \omega^{n, \theta_2} \, + \,  (\Theta ^n)^2 \ \mbox{for} \; \,  n=1, 2,\dots N,
			\end{split}
		\end{equation*}
		where $\omega^{k, \theta_i} := \theta_i \omega^{k-1} + (1- \theta_i) \omega^{k}$ for $\theta_i \in [0,1],\  i = 1, 2$.  Then
		\begin{equation*}
			\begin{split}
				\omega^{n} \le 2 E_{\alpha} \bigg( \frac{11\Lambda \, t_{n}^{\alpha}}{2}\bigg) \, &\Big[ \omega^0 + \max_{1 \le k \le n} \sum_{j=1}^{k} \, Q^{(k)}_{k-j} (\zeta^j + \Theta^j) + \max_{1 \le j \le n} \, \left\lbrace \Theta^j \right\rbrace \Big].\\
			\end{split}
		\end{equation*}
		Provided the maximum time-step $\tau_N \le \frac{1}{\sqrt[\alpha]{\frac{11}{2} \, \Gamma (2- \alpha) \, \Lambda}}.$
	\end{lemma}

\begin{theorem}\label{STL1}
		The solution $(U_h^n , V_h^n)$ of the fully-discrete scheme \eqref{fd}  satisfies the following estimates
		\begin{equation*}
				\begin{split}
					\max_{1 \le n \le N}\|U_h^{n}\| &\le C \, \big(\|U_h^{0}\| \big),\\
					\max_{1 \le n \le N}\|V_h^{n}\| &\le C \, \big(\|V_h^{0}\|+\|U_h^{0}\| \big).
				\end{split}
		\end{equation*}
	\end{theorem}	
	
\noindent\textbf{Proof.}  First we consider the solution $U_h^n$ of the proposed scheme for $n\geq 2.$ Observe that the equation \eqref{fd2} provides
	\begin{equation}\label{St1}
		( V_h^{n,\sigma}, z_h)-(\nabla U_h^{n,\sigma},\nabla z_h)=0 \quad \mbox{for} \ n \geq 2.
	\end{equation}
	On adding equations \eqref{fd1} and \eqref{St1} with $y_h=U_h^{n,\sigma}$ and $z_h=V_h^{n,\sigma}$, we obtain
	\begin{equation}\label{stt3}
		(({D_N^\alpha}U_h)^{n-\sigma},U_h^{n,\sigma})+\|V_h^{n,\sigma}\|^2+( V_h^{n,\sigma}, U_h^{n,\sigma})
		=(f(\hat{U}_h^{n-1,\sigma} ),U_h^{n,\sigma}). 
	\end{equation}
	Again utilizing \eqref{St1} with $z_h=U_h^{n,\sigma},$ \eqref{lip},  and Poincar$\acute{e}$ inequality in \eqref{stt3}   to achieve 
	\begin{equation}\label{stt4}
 \begin{split}
		(({D_N^\alpha}U_h)^{n-\sigma},U_h^{n,\sigma})+\|V_h^{n,\sigma}\|^2+\|\nabla U_h^{n,\sigma}\|^2
		&=(f(\hat{U}_h^{n-1,\sigma} ),U_h^{n,\sigma})\\
            &\le L_{f}\|\hat{U}_h^{n-1,\sigma}\| \|U_h^{n,\sigma}\|\\
            &\le \frac{(L_{f}C_p)^2}{4}\|\hat{U}_h^{n-1,\sigma}\|^2+ \|\nabla U_h^{n,\sigma}\|^2,\\
            ({D_N^\alpha}\|U_h\|^2)^{n-\sigma}
			\le (L_{f}C_p(1+C_r))^2&\|{U}_h^{n-1}\|^2+(C_r L_{f}C_p)^2\|{U}_h^{n-2}\|^2 \ \text{for}\  n\geq 2,
   \end{split}
	\end{equation}
	where we have used Lemma \ref{l111} and the following facts
 \begin{equation}\label{fact}
     \begin{split}
		\tau_{n,\sigma}^*&\le (1-\sigma)\frac{\tau_n}{\tau_{n-1}}\le\frac{\tau_n}{\tau_{n-1}}\le C_r,\\
				\|\hat{U}_h^{n-1,\sigma}\|^2&\le 2(1+\tau_{n,\sigma}^*)^2\|{U}_h^{n-1}\|^2+2(\tau_{n,\sigma}^*)^2\|{U}_h^{n-2}\|^2\\
    &\le 2(1+C_r)^2\|{U}_h^{n-1}\|^2+2C_r^2\|{U}_h^{n-2}\|^2.
			\end{split}
	\end{equation}
	
 
\noindent Next, we consider the solution $U_h^n$ of the proposed scheme for $n=1.$ Observe that the equation
 \eqref{fd2} provides
\begin{equation}\label{stt18}
		( V_h^{1,\sigma}, z_h)-(\nabla U_h^{1,\sigma},\nabla z_h)=0.
	\end{equation}
	Following the similar line proof as above of the solution  $U_h^n$ for $n\geq 2$ after adding equations \eqref{fnd1} and \eqref{stt18} with $x_h=U_h^{1,\sigma}$ and $z_h=V_h^{1,\sigma}$ yields
	\begin{equation}\label{stt24}
		({D_N^\alpha}\|U_h\|^2)^{1-\sigma}
		\le 4L_{f}(1-\sigma)^2\|{U}_h^{1}\|^2+4L_f\sigma^2\|{U}_h^{0}\|^2. 
	\end{equation}
	

\noindent	Here, first we combine inequities \eqref{stt24} and \eqref{stt4} for each $n\geq 1$ and then apply the Lemma \ref{l5}, we get the first required result
	\begin{equation}\label{stt28}
		\max_{1 \le n \le N}\|U_h^{n}\| \le C\,\|U_h^0\|.
	\end{equation}
	
	\noindent Now, we prove the next result for the solution $V_h^n \ \  \forall\,n\geq 1.$
 Note that from presented scheme \eqref{fd2} one has the following equation for each $j=1,2,\dots,n$
	\begin{equation*}
		(V_h^j-V_h^{j-1},z_h)-(\nabla U_h^j-\nabla U_h^{j-1}, \nabla z_h)=0.
	\end{equation*}
	Further multiplying by $\mathcal{A}_{\sigma,n-j}^{n}$ on both side and summing over $j=1,2,\dots,n$ to get 
	\begin{equation}\label{stt30}
		((D_N^\alpha V_h)^{n-\sigma},z_h)-((D_N^\alpha \nabla U_h)^{n-\sigma}, \nabla z_h)=0.
	\end{equation}
On adding equations \eqref{fd1} and \eqref{stt30} with  $y_h=-\Delta_h V_h^{n,\sigma}$ and  $z_h=V_h^{n,\sigma}$, we obtain	
	\begin{equation}\label{stt34}
		(({D_N^\alpha}V_h)^{n-\sigma}, V_h^{n,\sigma}) \le\frac{L_{f}}{4}
		 \|\hat{U}_h^{n-1,\sigma}\|^2. 
	\end{equation}
 \noindent Using Lemma \ref{l111} and \eqref{stt28} in \eqref{stt34}, we reach at 
	\begin{equation}\label{stt35}
		(\|D_N^\alpha \|V_h\|^2)^{n-\sigma} \le C\,\|U_h^0\| \ \mbox{for} \ n\geq2.
	\end{equation}
 
	\noindent  Next, we consider the solution $V_h^n$ of the proposed scheme for $n=1.$ Observe that the equation \eqref{stt30} provides
	\begin{equation}\label{stt39}
		((D_N^\alpha V_h)^{1-\sigma},z_h)-((D_N^\alpha \nabla U_h)^{1-\sigma}, \nabla z_h)=0.
	\end{equation}
Following the similar line proof as above of the solution  $V_h^n$ for $n\geq 2$ after adding equations \eqref{fnd1} and \eqref{stt39} with $x_h=-\Delta_h V_h^{1,\sigma}$ and $z_h=V_h^{1,\sigma}$  yields	
\begin{equation}\label{stt44}
(D_N^\alpha \|V_h\|^2)^{1-\sigma}
		\le C\|U_h^0\|^2. 
\end{equation}
Finally, on combining inequalities \eqref{stt44} and \eqref{stt35} for each $n\geq 1$ and then apply the Lemma \ref{l5}, we get the next required result
 \begin{equation}\label{stt47}
		\max_{1 \le n \le N}\|V_h^n\| \le C\big[\|V_h^0\|+\|U_h^0\|\big].
	\end{equation}
	\hfill{$\Box$}

\noindent Furthermore, to obtain error estimate for the fully-discrete solution $(U_h^n, V_h^n)$, we require certain reasonable regularity assumptions on the solution $u$ as well as some more notations and results that are given below:
	\begin{equation}\label{res}
		\begin{split}
			\|{D}^{\alpha}_{t}u\|_{L^{\infty}(0,T; H_0^{1}(\Omega)\cap H^{m+3}(\Omega))}\le C \ \   \mbox{and} \ &\|\partial_t^q u(\cdot,t)\|_{H^{m+3}(\Omega)}\le  C( 1 + t^{\alpha - q}),\\
   \ &\text{for} \ t\in (0,T],\  \mbox{where} \ q=0,\,1,\,2,\,\mbox{and} \  3.
		\end{split}
	\end{equation}
	The following property for the Ritz projection $R_h$ holds in more general spaces.
	\begin{theorem}\label{proj}\cite{[vth],[JC]}
		There exists a constant $C$ that is positive and independent of $h$ such that
		\begin{equation*}
			||\chi(t)||_{H^l(\Omega)}\leq  C(w,\,\Omega,\,m)h^{k-l}||u(t)||_{H^k(\Omega)} \quad\quad  \forall \,  u(t)\in H^k(\Omega)\cap H_0^1(\Omega),
		\end{equation*}
        \qquad $\text{for}\ l=0,1 \ \text{and} \ l<k\leq m+1,$
		where $\chi(t)=u(t)-R_h u(t).$
	\end{theorem}
	\begin{lemma}\label{cl6}
		\cite{[hr12]} If $\|\partial_t^q \omega(\cdot,t)\|_{H^2(\Omega)}\le C( 1 + t^{\alpha-q })$ for $q=0,1,2,3$ and $t \in (0,T].$ Then for every  $n=1, 2,\dots, N$, the local truncation error satisfies
		\begin{equation*}
			\, \| D_{t_{n-\sigma}}^{\alpha}\omega- D^{\alpha}_{N} \, \omega^{n- \sigma} \|_{H^2(\Omega)} \le  Ct_{n-\sigma}^{-\alpha} \, N^{-\min \{3- \alpha, \, r \alpha\}}.
		\end{equation*}
	\end{lemma}  
	\begin{lemma}\label{l3}
		\cite{[hr12]} If $\|\partial_t^q \omega(\cdot,t)\|_{H^2(\Omega)} \le C \, (1+t^{\alpha-q})$ for $q = 0,1,2,$   for each  $n=1, 2,\dots, N,$ we have
		\begin{equation*}
				 \|\omega^{n, \sigma} -\omega^{n-\sigma} \|_{H^2(\Omega)} \le  C N^{-\min\{r \alpha, \, 2\}}  \ \mbox{and} \  \|\hat{\omega}^{{n-1}, \sigma} -\omega^{n-\sigma} \| \le  C N^{-\min\{r \alpha, \, 2\}}. 
     		\end{equation*}
	\end{lemma}
 
	\noindent Now, we prove the error estimate for the fully-discrete solution using the $\alpha$-robust Gronwall inequality given in Lemma \ref{l5}. For this purpose, we rewrite the errors $u^n-U_h^{n}$ and $v^n-V_h^{n}$ using the Ritz projection $R_h$ as
	\begin{equation*}\label{cufc}
		\begin{split}
			u^n-U_h^n=&(u^n-R_hu^n)+(R_hu^n-U_h^n)\\
			=& \chi_u^n +\eta_u^n.\\
			v^n-V_h^n=&(v^n-R_hv^n)+(R_hv^n-V_h^n)\\
			=& \chi_v^n +\eta_v^n.
		\end{split}
	\end{equation*}
 The exact solution $(u^n,v^n)$ of problem \eqref{bha1} satisfies the following equation with error terms  $E_c^{n-\sigma}$  and  $E_f^{n-\sigma}$ corresponding to Caputo derivative approximation and nonlinear source term $f$ respectively.
 \begin{equation}\label{ebha2}
 \begin{split}
          (D_N^\alpha u)^{n-\sigma}-\Delta v^{n-\sigma}+v^{n-\sigma}&= \  f(\hat{u}^{n-1,\sigma})+E_c^{n-\sigma}+E_f^{n-\sigma}\ \ \mbox{for} \  n \geq 2,\\
       (D_N^\alpha u)^{1-\sigma}-\Delta v^{1-\sigma}+v^{1-\sigma}&= \  f({u}^{1,\sigma})+E_c^{1-\sigma}+E_f^{1-\sigma} \ \ \ \ \   \mbox{for} \  n=1, 
 \end{split}
 \end{equation}
 
\noindent where  $E_c^{n-\sigma}= (D_N^\alpha u)^{n-\sigma}-D_{t_{n-\sigma}}^\alpha u \  \  \mbox{if} \  n\geq 1,$ and
\begin{equation*}
 E_f^{n-\sigma}= \begin{cases} f(u^{n-\sigma})-f(\hat{u}^{n-1,\sigma}) & \mbox{} \  \  \mbox{if} \  n\geq2, \\ \,  f(u^{1-\sigma})-f({u}^{1,\sigma})  & \mbox{} \  \  \mbox{if} \ n=1.\\\end{cases}
 \end{equation*}
 \noindent Use \eqref{lip} and Lemma \ref{l3}  to get
	\begin{equation}\label{Ef}
 \begin{split}
    \|E_f^{n-\sigma}\|\le C\, N^{-\min \{2, \, r \alpha\}} \ \mbox{for} \ n\geq1.
 \end{split}
	\end{equation}
	\begin{theorem}\label{th8}
		The solution $(u^n,v^n)$ of \eqref{bha1} and the solution $(U_h^n , V_h^n)$ of fully-discrete scheme \eqref{fd} satisfies the following inequalities for each  $n=1, 2,\dots, N$ 
		\begin{equation*}
		    \begin{split}
		        \max_{1 \le n \le N} \|u^{n} - U_h^{n}\| &\le  C \big(h^{m+1} + N^{-\min \left\lbrace 2, \, r \alpha \right\rbrace }\big), \\
          	\max_{1 \le n \le N} \|v^{n} - V_h^{n}\| &\le  C \big(h^{m+1} + N^{-\min \left\lbrace 2, \, r \alpha \right\rbrace }\big). 
		    \end{split}
		\end{equation*}
	\end{theorem}
	\noindent \textbf{Proof.} 
	 First we consider the following equation  for $n\geq2$
	\begin{equation*}
		\begin{split}
			&\left( (D_N^\alpha \eta_u)^{n-\sigma}, y_h \right) \, + \, (\nabla \eta_v^{n,\sigma}, \nabla y_h)+(\eta_v^{n,\sigma}, y_h)  \qquad\\
   &\;=\left( (D_N^\alpha (R_hu-U_h))^{n-\sigma}, y_h \right) \, + \, (\nabla (R_h v^{n,\sigma}-V_h^{n,\sigma}), \nabla y_h)+(R_hv^{n,\sigma}-V_h^{n,\sigma}, y_h)\\
			&\;=\left( (D_N^\alpha R_hu)^{n-\sigma}, y_h \right) \, + \, (\nabla v^{n,\sigma}, \nabla y_h)+(R_hv^{n,\sigma}, y_h)-(f(\hat{U}_h^{n-1,\sigma}),y_h).
		\end{split}
	\end{equation*}
 Application of equation \eqref{ebha2} yields 
 \begin{equation}\label{e2}
		\begin{split}
  &\left( (D_N^\alpha \eta_u)^{n-\sigma}, y_h \right) \, + \, (\nabla \eta_v^{n,\sigma}, \nabla y_h)+(\eta_v^{n,\sigma}, y_h)\\
   &\;=\left( (D_N^\alpha (R_hu-u))^{n-\sigma}, y_h \right) \, - \, (\Delta (v^{n,\sigma}- v^{n-\sigma}),  y_h)+(R_hv^{n,\sigma}-v^{n-\sigma}, y_h)\\
			&\;\qquad +(E_c^{n-\sigma},y_h)+(E_f^{n-\sigma},y_h)+\Big(f(\hat{u}^{n-1,\sigma})-f(\hat{U}_h^{n-1,\sigma}),y_h\Big)\\
			&\;=\left(- (D_N^\alpha \chi_u)^{n-\sigma}, y_h \right) \, - \, (\Delta (v^{n,\sigma}- v^{n-\sigma}),  y_h)+(v^{n,\sigma}-v^{n-\sigma}-\chi_v^{n,\sigma}, y_h)\\
			&\;\qquad +(E_c^{n-\sigma},y_h)+(E_f^{n-\sigma},y_h)+\Big(f(\hat{u}^{n-1,\sigma})-f(\hat{U}_h^{n-1,\sigma}),y_h\Big).
		\end{split}
	\end{equation}
	\noindent	Next, we consider the following equation for $n\geq1$	
	\begin{equation}\label{e3}
		\begin{split}
			&(\eta_v^n, z_h)-(\nabla \eta_u^n, \nabla z_h)=(R_hv^{n},z_h)-(\nabla R_hu^n, \nabla z_h)-(V_h^{n},z_h)+(\nabla U_h^n, \nabla z_h)\\
			&\;\qquad\qquad\qquad\qquad\quad=(R_hv^{n},z_h)-(\nabla u^n, \nabla z_h) \\
			&\;\qquad\qquad\qquad\qquad\quad=(- \chi_v^n, z_h).
		\end{split}
	\end{equation}
	\noindent Further, the equation \eqref{e3} leads to  
	\begin{equation}\label{e4}
		(\eta_v^{n,\sigma}, z_h)-(\nabla \eta_u^{n,\sigma}, \nabla z_h)=(- \chi_v^{n,\sigma}, z_h). 
	\end{equation}

	\noindent On adding equations \eqref{e2} and \eqref{e4} with $y_h=\eta_u^{n,\sigma}$ and $z_h=\eta_v^{n,\sigma}$, we obtain 
	\begin{equation}\label{e7}
		\begin{split}
			&( (D_N^\alpha \eta_u)^{n-\sigma},\eta_u^{n,\sigma} ) \, +(\eta_v^{n,\sigma}, \eta_u^{n,\sigma})+\|\eta_v^{n,\sigma}\|^2\\
			&=(- (D_N^\alpha \chi_u)^{n-\sigma},\eta_u^{n,\sigma} ) \, - \, (\Delta (v^{n,\sigma}- v^{n-\sigma}),  \eta_u^{n,\sigma})+(v^{n,\sigma}-v^{n-\sigma}-\chi_v^{n,\sigma}, \eta_u^{n,\sigma})\\
			&\;\quad+(E_c^{n-\sigma},\eta_u^{n,\sigma})+(E_f^{n-\sigma},\eta_u^{n,\sigma})+\Big(f(\hat{u}^{n-1,\sigma})-f(\hat{U}_h^{n-1,\sigma}),\eta_u^{n,\sigma}\Big)+(- \chi_v^{n,\sigma}, \eta_v^{n,\sigma}).
		\end{split}
	\end{equation}
	
	\noindent Again utilizing \eqref{e4} with $z_h=\eta_u^{n,\sigma}$ and \eqref{lip} in \eqref{e7} to reach at
	\begin{equation*}
		\begin{split}
			( (D_N^\alpha \eta_u)^{n-\sigma},&\eta_u^{n,\sigma} ) \, +\|\nabla \eta_u^{n,\sigma}\|^2+\|\eta_v^{n,\sigma}\|^2\\
   	&\le\big[\|(D_N^\alpha \chi_u)^{n-\sigma}\|  \, + \, \|\Delta (v^{n,\sigma}- v^{n-\sigma})\|+\|v^{n,\sigma}-v^{n-\sigma}\|+L_f\|\hat{\chi_u}^{n-1,\sigma}\| \\
			&\qquad+\|E_c^{n-\sigma}\|+\|E_f^{n-\sigma}\|\big]\|\eta_u^{n,\sigma}\|+L_f\|\hat{\eta_u}^{n-1,\sigma}\|\|\eta_u^{n,\sigma}\|+ \|\chi_v^{n,\sigma}\|\|\eta_v^{n,\sigma}\|\\
   &\le\big[\|(D_N^\alpha \chi_u)^{n-\sigma}\|  \, + \, \|\Delta (v^{n,\sigma}- v^{n-\sigma})\|+\|v^{n,\sigma}-v^{n-\sigma}\|+L_f\|\hat{\chi_u}^{n-1,\sigma}\| \\
			&\qquad+\|E_c^{n-\sigma}\|+\|E_f^{n-\sigma}\|\big]\|\eta_u^{n,\sigma}\|+\frac{(L_fC_p)^2}{4}\|\hat{\eta_u}^{n-1,\sigma}\|^2+\|\nabla \eta_u^{n,\sigma}\|^2\\
   &\hspace{2cm}+ \frac{1}{4}\|\chi_v^{n,\sigma}\|^2+\|\eta_v^{n,\sigma}\|^2.
		\end{split}
	\end{equation*}
 Further, the Lemma \ref{l111} and the facts given in \eqref{fact} produce
 \begin{equation*}
     \begin{split}
         (D_N^\alpha \|\eta_u\|^2)^{n-\sigma} &\le 2\big[\|(D_N^\alpha \chi_u)^{n-\sigma}\|  \, + \, \|\Delta (v^{n,\sigma}- v^{n-\sigma})\|+\|v^{n,\sigma}-v^{n-\sigma}\|\\
   &\qquad+L_f\|\hat{\chi_u}^{n-1,\sigma}\|+\|E_c^{n-\sigma}\|+\|E_f^{n-\sigma}\|\big]\|\eta_u^{n,\sigma}\|\\
&\qquad+(L_{f}C_p(1+C_r))^2\|\eta_u^{n-1}\|^2+(L_{f}C_pC_r)^2\|\eta_u^{n-2}\|^2+ \frac{1}{2}\|\chi_v^{n,\sigma}\|^2.
     \end{split}
 \end{equation*}
 By making use of \eqref{Ef} together with  Lemmas \ref{cl6} and \ref{l3}, we have 
 \begin{equation*}
     \begin{split}
         (D_N^\alpha \|\eta_u\|^2)^{n-\sigma} &\le 2\big[\|(D_N^\alpha \chi_u)^{n-\sigma}\|  +C\left(N^{-\min \{2, \, r \alpha\}}+t_{n-\sigma}^{-\alpha} \, N^{-\min \{3- \alpha, \, r \alpha\}}\right)\\
   &\qquad+L_f\|\hat{\chi_u}^{n-1,\sigma}\|\Big]\|\eta_u^{n,\sigma}\|+(L_{f}C_p(1+C_r))^2\|\eta_u^{n-1}\|^2\\
&\qquad+(L_{f}C_pC_r)^2\|\eta_u^{n-2}\|^2+\frac{1}{2}\|\chi_v^{n,\sigma}\|^2.
     \end{split}
 \end{equation*}
 Finally, Theorem \ref{proj} gives the following inequality for each $n\geq2$
\begin{equation}\label{e11}
		\begin{split}
			&(D_N^\alpha \|\eta_u\|^2)^{n-\sigma} \le C\Big[\Big(h^{m+1}+N^{-\min \{2, \, r \alpha\}}+t_{n-\sigma}^{-\alpha} \, N^{-\min \{3- \alpha, \, r \alpha\}}\Big)\|\eta_u^{n,\sigma}\|\\
   &\qquad\qquad\qquad\qquad\qquad\qquad\qquad\qquad\quad+\|{\eta}_u^{n-1}\|^2
       + \|{\eta}_u^{n-2}\|^2+h^{m+1}\Big].
		\end{split}
	\end{equation}
\noindent Similar to above, we can write the following equations for $n=1$
 \begin{equation}\label{u1e}
		\begin{split}
			&\left( (D_N^\alpha \eta_u)^{1-\sigma}, x_h \right) \, + \, (\nabla \eta_v^{1,\sigma}, \nabla x_h)+(\eta_v^{1,\sigma}, x_h)  \qquad\\
             &\;=\left(- (D_N^\alpha \chi_u)^{1-\sigma}, x_h \right) \, - \, (\Delta (v^{1,\sigma}- v^{1-\sigma}),  x_h)+(v^{1,\sigma}-v^{1-\sigma}-\chi_v^{1,\sigma}, x_h)\\
			&\;\qquad\qquad\qquad+(E_c^{1-\sigma}, x_h)+(E_f^{1-\sigma}, x_h)+\Big(f(u^{1,\sigma})-f(U_h^{1,\sigma}),x_h\Big),
		\end{split}
	\end{equation}
	\begin{equation*}
		(\eta_v^{1,\sigma}, z_h)-(\nabla \eta_u^{1,\sigma}, \nabla z_h)=(- \chi_v^{1,\sigma}, z_h). 
	\end{equation*}

\noindent Now, applying the parallel argument as in the case of $n\geq2$, we obtain 
\begin{equation}\label{e11n}
\begin{split}
   &(D_N^\alpha \|\eta_u\|^2)^{1-\sigma} \le C\Big[\Big(h^{m+1}+N^{-\min \{2, \, r \alpha\}}+t_{1-\sigma}^{-\alpha} \, N^{-\min \{3- \alpha, \, r \alpha\}}\Big)\|\eta_u^{1,\sigma}\|\\
  &\qquad\qquad\qquad\qquad\qquad\qquad\qquad\qquad\qquad+\|\eta_u^{1}\|^2+\|\eta_u^{0}\|^2+h^{m+1}\Big],
\end{split}
\end{equation}

	\noindent On combining \eqref{e11n} and \eqref{e11} for each $n\geq 1$ and then apply the Lemma \ref{l5} provides
 \begin{equation*}
		\begin{split}
			\|\eta_u^{n}\| \le& C \, \Big[ \|\eta_u^{0}\| + \max_{1 \le k \le n} \sum_{j=1}^{k} \, Q^{(k)}_{k-j} \left(h^{m+1}+N^{-\min \{2, \, r \alpha\}}\right)+\\
            &\qquad\qquad\max_{1 \le k \le n} \sum_{j=1}^{k} \, Q^{(k)}_{k-j}\left(t_{j-\sigma}^{-\alpha} \, N^{-\min \{3- \alpha, \, r \alpha\}}\right) +h^{m+1}\Big].
		\end{split}
	\end{equation*}
	
	\noindent Notice that, if $l_N=\frac{1}{\ln{N}}$ then for each $j=1,\,2,\,\dots,n$ one has 
	\begin{equation}\label{e127}
		t_{j-\sigma}^{-\alpha}\le  (j/2)^{-r\alpha}N^{r\alpha}\le(2N)^{r\alpha}j^{r(l_N-\alpha)}.
	\end{equation}
   Now, we use Lemma \ref{l7} along with  $\gamma=\alpha$ and  \eqref{e127} to acquire
	\begin{equation}\label{u1ee14}
		\|\eta_u^{n}\| \le C\big(\|\eta_u^{0}\|+h^{m+1}+N^{-min\{2,r\alpha\}}\big) \ \ \mbox{for} \ \  n\geq 1.
	\end{equation}

	\noindent At last, triangle inequality and Theorem \ref{proj} together with $U_h^0=R_h u_0$ gives the required result  
	\begin{equation*}
		\max_{1 \le n \le N} \|u^{n} - U_h^{n}\| \le  C \big(h^{m+1} + N^{-\min \left\lbrace 2, \, r \alpha \right\rbrace }\big). 
	\end{equation*}
 Now, we derive the error estimate for the fully-discrete solution $V_h^n$. Further, the following equation can be obtained from \eqref{e3} by using same approach as in \eqref{stt30} 
	\begin{equation}\label{e15}
		((D_N^\alpha \eta_v)^{n-\sigma}, z_h)-((D_N^\alpha \nabla \eta_u)^{n-\sigma}, \nabla z_h)=(- (D_N^\alpha\chi_v)^{n-\sigma}, z_h) \ \mbox{for} \ n\geq1. 
	\end{equation}
\noindent On adding equations \eqref{e2} and \eqref{e15} with $y_h=-\Delta_h \eta_v^{n,\sigma}$ and $z_h=\eta_v^{n,\sigma}$, we obtain
\begin{equation}\label{e16d}
\begin{split}
	( (D_N^\alpha \nabla &\eta_u )^{n-\sigma},\nabla \eta_v^{n,\sigma} ) \, +\|\Delta_h \eta_v^{n,\sigma}\|^2+\|\nabla \eta_v^{n,\sigma}\|^2\\
     &=\left(- (D_N^\alpha \chi_u)^{n-\sigma}, -\Delta_h \eta_v^{n,\sigma} \right) \, - \, (\Delta (v^{n,\sigma}- v^{n-\sigma}),  -\Delta_h \eta_v^{n,\sigma})\\
   &\quad +(v^{n,\sigma}-v^{n-\sigma}-\chi_v^{n,\sigma}, -\Delta_h \eta_v^{n,\sigma})+(E_c^{n-\sigma}, -\Delta_h \eta_v^{n,\sigma})+(E_f^{n-\sigma}, -\Delta_h \eta_v^{n,\sigma})\\
   &\quad+\Big(f(\hat u^{n-1,\sigma})-f(\hat U_h^{n-1,\sigma}),-\Delta_h \eta_v^{n,\sigma}\Big)+( (D_N^\alpha\chi_v)^{n-\sigma},\eta_v^{n,\sigma})\\
   &\ \le \left(- (D_N^\alpha \chi_u)^{n-\sigma}, -\Delta_h \eta_v^{n,\sigma} \right)-(E_c^{n-\sigma}, \Delta_h \eta_v^{n,\sigma}) +\frac{1}{4}\Big[ \|\Delta (v^{n,\sigma}- v^{n-\sigma})\|\\
&\quad +\|v^{n,\sigma}-v^{n-\sigma}\|+\|\chi_v^{n,\sigma}\|+\left\|f\left(\hat{u}^{n-1,\sigma}\right)-f\left(\hat{U_h}^{n-1,\sigma}\right)\right\|+\|E_f^{n-\sigma}\|\Big]^{2}\\
&\quad+\|\Delta_h \eta_v^{n,\sigma}\|^2+\frac{C_p^2}{4}\|(D_N^\alpha \chi_v)^{n-\sigma}\|^2+\|\nabla \eta_v^{n,\sigma}\|^2.
   \end{split}
\end{equation} 
Note that 
\begin{equation}\label{fact2}
    \left(- (D_N^\alpha \chi_u)^{n-\sigma}, -\Delta_h \eta_v^{n,\sigma} \right)=\left( E_c^{n-\sigma}, \Delta_h \eta_v^{n,\sigma} \right)-\left(R_h E_c^{n-\sigma}, \Delta_h \eta_v^{n,\sigma} \right)+\left( D_{t_{n-\sigma}}^\alpha \chi_u,\Delta_h \eta_v^{n,\sigma} \right)\\
\end{equation}
Further, the use of \eqref{pr} and \eqref{DK25} in \eqref{fact2} provide
\begin{equation}\label{fact3}
   \left(- (D_N^\alpha \chi_u)^{n-\sigma}, -\Delta_h \eta_v^{n,\sigma} \right)=\left( E_c^{n-\sigma}, \Delta_h \eta_v^{n,\sigma} \right)-\left(\Delta  E_c^{n-\sigma},  \eta_v^{n,\sigma} \right)+\left( D_{t_{n-\sigma}}^\alpha \chi_u,\Delta_h \eta_v^{n,\sigma} \right)
\end{equation}
Through the application of \eqref{fact3} in \eqref{e16d}, it can be deduced that 
\begin{equation*}
\begin{split}
 (D_N^\alpha \|\eta_v\|^2)^{n-\sigma}&\le2\|\Delta E_c^{n-\sigma}\|\| \eta_v^{n,\sigma}\|+\frac{1}{2}\big[\|D_{t_{n-\sigma}}^\alpha \chi_u\|  \, + \, \|\Delta (v^{n,\sigma}- v^{n-\sigma})\|+\|v^{n,\sigma}-v^{n-\sigma}\|\\
&\;+\|\chi_v^{n,\sigma}\|+L_f\|\hat{\chi_u}^{n-1,\sigma}\|+L_f\|\hat{\eta_u}^{n-1,\sigma}\| +\|E_f^{n-\sigma}\|+C_p\|(D_N^\alpha \chi_v)^{n-\sigma}\|\big]^{2}.
   \end{split}
\end{equation*} 
 By making use of \eqref{Ef} along with  Lemmas \ref{cl6} and \ref{l3}, we obtain
 \begin{equation*}
\begin{split}
 (D_N^\alpha \|\eta_v\|^2)^{n-\sigma}&\le C\big[t_{n-\sigma}^{-\alpha} \, N^{-\min \{3- \alpha, \, r \alpha\}} \| \eta_v^{n,\sigma}\|+\big(\|D_{t_{n-\sigma}}^\alpha \chi_u\|  \, + N^{-\min \{2, \, r \alpha\}}\\
&\;+\|\chi_v^{n,\sigma}\|+\|\hat{\chi_u}^{n-1,\sigma}\|+\|\hat{\eta_u}^{n-1,\sigma}\| +\|(D_N^\alpha \chi_v)^{n-\sigma}\|\big)^{2}\big].
   \end{split}
\end{equation*} 
 Finally, Theorem \ref{proj} and \eqref{u1ee14} give the following inequality for each $n\geq2$
\begin{equation}\label{e19n}
		\begin{split}
			& (D_N^\alpha \|\eta_v\|^2)^{n-\sigma}\le C\Big[t_{n-\sigma}^{-\alpha} \, \big(N^{-\min \{3- \alpha, \, r \alpha\}} \big)\|\eta_v^{n,\sigma}\|+\big[h^{m+1}+\big(N^{-\min \{2, \, r \alpha\}}\big)\big]^{2}\Big].
		\end{split}
	\end{equation}

\noindent Observe that from the equation \eqref{e15}, we have
	\begin{equation}\label{v1e15}
		((D_N^\alpha \eta_v)^{1-\sigma}, z_h)-((D_N^\alpha \nabla \eta_u)^{1-\sigma}, \nabla z_h)=(- (D_N^\alpha\chi_v)^{1-\sigma}, z_h). 
	\end{equation}
 Following the similar line proof as above in the case of $n\geq2$ after adding equations \eqref{u1e} and \eqref{v1e15} with $x_h=-\Delta_h \eta_v^{1,\sigma}$ and $z_h=\eta_v^{1,\sigma}$  yields	

\begin{equation}\label{v1e19}
		(D_N^\alpha \|\eta_v\|^2)^{1-\sigma}\le C\Big[t_{1-\sigma}^{-\alpha} \, \big(N^{-\min \{3- \alpha, \, r \alpha\}} \big)\|\eta_v^{1,\sigma}\|+\big[h^{m+1}+\big(N^{-\min \{2, \, r \alpha\}}\big)\big]^{2}\Big].
	\end{equation} 

 \noindent  On combining \eqref{e19n} and \eqref{v1e19} for each $n\geq 1$ and then apply Lemma \ref{l5} provides
	\begin{equation*}
		\begin{split}
			\|\eta_v^{n}\| \le& C\, \Big[ \|\eta_v^{0}\| +  \max_{1 \le k \le n} \sum_{j=1}^{k} \, Q^{(k)}_{k-j}  \left(t_{j-\sigma}^{-\alpha} \, N^{-\min \{3- \alpha, \, r \alpha\}} \right)\\
   &+\max_{1 \le k \le n} \sum_{j=1}^{k} \, Q^{(k)}_{k-j} \left(N^{-\min \{2, \, r \alpha\}} +h^{m+1}\right) +\Big(N^{-\min \{2, \, r \alpha\}} +h^{m+1}\Big) \Big].
		   \end{split}
	\end{equation*}
  Now, we use Lemma \ref{l7} along with  $\gamma=\alpha$ and  \eqref{e127} to reach at
	\begin{equation*}
		\|\eta_v^{n}\| \le C\big(\|\eta_v^{0}\|+h^{m+1}+N^{-min\{2,r\alpha\}}\big) \ \ \mbox{for} \ \  n\geq 1.
	\end{equation*}
\noindent At last, triangle inequality and Theorem \ref{proj} together with $V_h^0=R_h v_0$ gives the next required result 
	\begin{equation*}
		\max_{1 \le n \le N} \|v^{n} - V_h^{n}\| \le  C \big(h^{m+1} + N^{-\min \left\lbrace 2, \, r \alpha \right\rbrace }\big). 
	\end{equation*}
	\hfill{$\Box$}

			\section{Numerical Experiments}
\begin{figure}[H]
	\begin{center}
		\includegraphics[width=13cm]{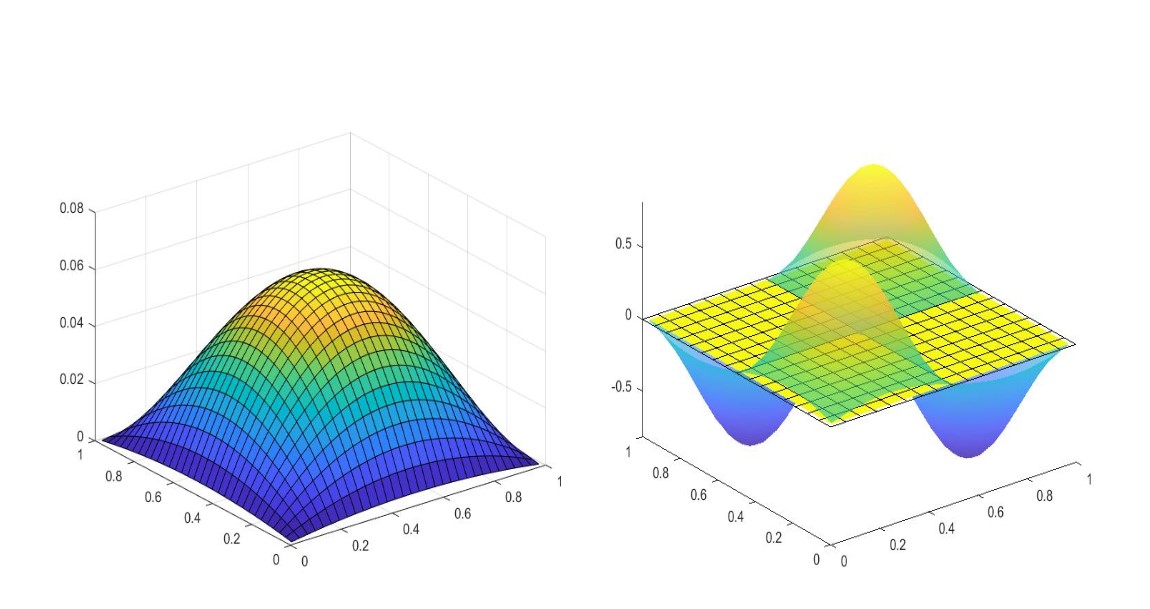}
	\end{center}
	\caption{\emph{Plot of the weight function (left) and  numerical solution with its domain (right) for Example 1. }}
	\label{pic01}
\end{figure}

	To validate the theoretical estimates, we perform two numerical experiments on different domains. We use scheme \eqref{fd} to solve time-fractional biharmonic equation with weakly singular solutions on two different domains. In both examples, we fix final time $T=0.5$ and set grading parameter $r=\frac{2}{\alpha}$ to achieve the optimal convergence rate in time. Initial guess $u^g$ is taken as the solution of the Poisson problem with the homogeneous Dirichlet boundary condition and tolerance is  $10^{-12}$ for Newton's iterations. To obtain spatial errors and corresponding convergence rates, we run the test for fixed $N=10000$ and $\alpha=0.5$. The temporal errors and convergence rates are calculated for $\alpha=0.4, \ 0.6,$ and $0.8$. Note that the established theoretical convergence rate for the fully-discrete scheme \eqref{fd} is $O\big
	(h^{m+1} + N^{-\min \left\lbrace 2, \, r \alpha \right\rbrace }\big)$ in  the $L^2(\Omega)$ norm.  Thus, temporal error is calculated by using relation $M= \left\lfloor{N^{\frac{2}{m+1}}}\right\rfloor$, where $M$ is the number of grid cells per coordinate direction. Now, denote the relative errors by \[E_u :=\frac{\displaystyle{\max_{1 \le n \le N}} \|u^{n} - U_h^{n}\|}{{\displaystyle{\max_{1\leq n\leq N}}\|u^{n}\|}} \qquad \mbox{and} \qquad E_v:=\frac{\displaystyle{\max_{1 \le n \le N}} \|v^{n} - V_h^{n}\|}{{\displaystyle{\max_{1\leq n\leq N}}\|v^{n}\|}}. \]\\ 
 {\textbf{Example 1}.}  In this example, we consider  domain $\Omega=(0,1)\times(0,1)$ and weight function $w=(1-x)(1-y)xy.$ We choose right hand side $f$ of \eqref{bha1} in such a way that solution is $u = (t^3+t^\alpha)(\sin(2\pi x)\sin(2\pi y))$. We also solve this example using  $L2$-$1_\sigma$ based standard finite element scheme to compare the effectiveness of the presented weighted $b$-spline based scheme \eqref{fd}. The results of the comparison between these two methods are presented in Table \ref{tab111} for weighted $b$-splines of degree $1$. Upon analyzing the errors and convergence rates in Table \ref{tab111}, it is clear that the weighted $b$-spline method performs better in terms of accuracy than the standard FEM. The numerical results obtained through the proposed $L2$-$1_{\sigma}$ method based scheme \eqref{fd} are presented in Tables \ref{tab5} and \ref{tab12}. In Table \ref{tab5}, we display errors and corresponding convergence rates for weighted $b$-splines of degrees $m=2$ and $3$. Further, Table \ref{tab12} reports temporal errors and their corresponding convergence rates. Finally, Figure \ref{pic01} illustrates the plot of the weight function and numerical solution at time $T=0.5 \ {\mbox{with}} \ M=16  \ \mbox{and} \ \alpha=0.5$.

	\begin{table}[H]
		\small
		\caption{Errors and corresponding rates in the spatial direction for Example 1.} \label{tab111}
		\begin{center}
			\begin{tabularx}{1\textwidth} { 
					>{\raggedright\arraybackslash}X 
					>{\centering\arraybackslash}X 
					>{\raggedleft\arraybackslash}X     >{\raggedleft\arraybackslash}X>{\raggedleft\arraybackslash}X>{\raggedleft\arraybackslash}X }
				\hline
				   & Degree of freedom & $E_u$  & Rate & $E_v$ & Rate  \\ 
				\hline
				Weighted  & 81 &  0.2257E-1  && 0.1714E-1 & \\  
				$b$-spline& 289 &0.5873E-2  & 1.9422&0.4353E-2 &1.9776  \\
			method	& 1089 &   0.1482E-2 & 1.9862& 0.1092E-2 &1.9950\\  
				error &4225  &   0.3714E-3 & 1.9967& 0.2732E-3 &1.9989\\	  
				\hline
				Standard & 81 &0.4255E-0 & &0.2217E-0 & \\  
			FEM	& 289 &0.1085E-0 &1.9714  & 0.5231E-1& 2.0834   \\
			 error	& 1089 & 0.2634E-1 &2.0423& 0.1246E-1 &2.0693 \\  
				& 4225 &0.6433E-2 &2.0338 & 0.3030E-2 &2.0402\\	  
				\hline
			\end{tabularx}
			
		\end{center}
	\end{table}
	
	\begin{table}[H]
		\small
		\caption{Errors  and corresponding rates in the spatial direction for Example 1.} \label{tab5}
		\begin{center}
			\begin{tabularx}{1\textwidth} { 
					>{\raggedright\arraybackslash}X 
					>{\centering\arraybackslash}X 
					>{\raggedleft\arraybackslash}X     >{\raggedleft\arraybackslash}X>{\raggedleft\arraybackslash}X>{\raggedleft\arraybackslash}X }
				\hline
				Degree   & Degree of freedom & $E_u$ & Rate & $E_v$& Rate  \\ 
				\hline 
				$m=2$ & 36 &   0.1278E-1  & & 0.1255E-1 & \\
				& 100 &   0.1423E-2  &3.1665& 0.1416E-2 &3.1478 \\  
				& 324 &0.1690E-3  & 3.0743 & 0.1688E-3 &3.0685   \\
				& 1156 &0.2084E-4 & 3.0193& 0.2084E-4 &3.0179\\  
				
				\hline
				$m=3$	& 49 &  0.2588E-2& & 0.2525E-2 & \\  
				& 121 &  0.1009E-3&4.6799 & 0.1009E-3 &4.6456 \\  
				& 361 &0.5511E-5 & 4.1956& 0.5511E-5 & 4.1944 \\
				& 1225 &  0.3321E-6 &4.0523& 0.3321E-6&4.0522\\  
				
				\hline
			\end{tabularx}
		\end{center}
	\end{table}
	
	\begin{table}[H]
		\small
		\caption{Errors and corresponding rates in the temporal direction for Example 1.} \label{tab12}
		\begin{center}
			\begin{tabularx}{1\textwidth} { 
					>{\raggedright\arraybackslash}X 
					>{\centering\arraybackslash}X 
					>{\raggedleft\arraybackslash}X     >{\raggedleft\arraybackslash}X>{\raggedleft\arraybackslash}X>{\raggedleft\arraybackslash}X  }
				\hline
				& $N$ & $E_{u}$  & Rate & $E_v$ & Rate  \\ 
				
				\hline
				$\alpha=0.4$	& 8 &  0.2964E-1  &  &0.2304E-1 & \\  
				& 16 &0.8707E-2 & 1.7675 &0.6856E-2 &1.7489   \\
				& 32 & 0.2343E-2 &1.8936&0.1867E-2 & 1.8762\\  
				&64 &0.6060E-3  &  1.9512 &  0.4862E-3 &1.9414\\	           
				\hline
				$\alpha=0.6$	& 8 &   0.2992E-1 &  &0.2331E-1 & \\  
				& 16 &0.8160E-2  & 1.8744 &0.6318E-2 &1.8836  \\
				& 32 &0.2104E-2 & 1.9550&0.1631E-2 & 1.9534\\  
				&64 & 0.5326E-3 & 1.9823 & 0.4134E-3 & 1.9801\\
				\hline
				$\alpha=0.8$	& 8 &   0.2970E-1  &  &0.2310E-1 & \\  
				& 16 &0.7842E-2  & 1.9213  &0.6007E-2 & 1.9433 \\
				& 32&0.1991E-2 & 1.9770&0.1521E-2 &1.9815 \\  
				&64 &0.5005E-3  & 1.9926 & 0.3820E-3 & 1.9933\\
				\hline
			\end{tabularx}
		\end{center}
	\end{table}

\begin{figure}[H]
	\begin{center}
		\includegraphics[width=13cm]{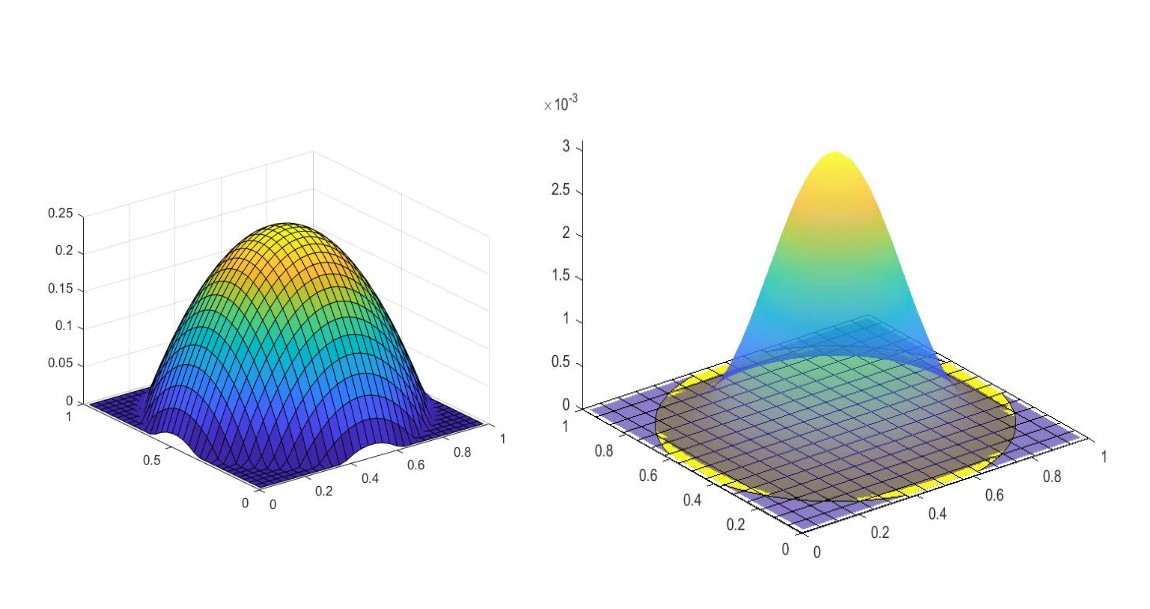}
	\end{center}
	\caption{\emph{Plot of the weight function (left) and  numerical solution with its domain (right) for Example 2. }}
	\label{pic02}
\end{figure}

	\noindent  {\textbf{Example 2}.}  In this example, we consider domain $\Omega=\big\{(x,y): \frac{1}{4}-\big(x-\frac{1}{2}\big)^2-\big(y-\frac{1}{2}\big)^2 > 0 \big\}$ and weight function $w=\frac{1}{4}-\big(x-\frac{1}{2}\big)^2-\big(y-\frac{1}{2}\big)^2.$ We choose right hand side $f$ of \eqref{bha1} in such a way that exact solution is $u = (t^3+t^\alpha)(w^3\sin(w))$.  The numerical results obtained through the proposed $L2$-$1_{\sigma}$ method based scheme \eqref{fd} are presented in Tables \ref{tab52} and \ref{tab122}. In Table \ref{tab52}, we show errors and corresponding convergence rates for weighted $b$-splines of degrees $m=1$, $2$, and $3$. Further,  Table \ref{tab122} reports temporal errors and their corresponding convergence rates. Figure \ref{pic02} illustrates the plot of the weight function and numerical solution at $T=0.5 \ \mbox{with} \ M=16 \ \mbox{and} \ \alpha=0.5$. 
\begin{table}[H]
		\small
		\caption{Errors and corresponding rates in the spatial direction for Example 2.} \label{tab52}
		\begin{center}
			\begin{tabularx}{1\textwidth} { 
					>{\raggedright\arraybackslash}X 
					>{\centering\arraybackslash}X 
					>{\raggedleft\arraybackslash}X     >{\raggedleft\arraybackslash}X>{\raggedleft\arraybackslash}X>{\raggedleft\arraybackslash}X }
				\hline
				Degree   & Degree of freedom & $E_u$ & Rate & $E_v$& Rate  \\ 
				\hline 
				$m=1$ & 81 &  0.5761E-1  & & 0.7484E-1 &\\
				& 289 &0.1536E-1  & 1.9069& 0.2004E-1 &1.9003   \\
				& 1089 &  0.3907E-2 &1.9753& 0.5104E-2 &1.9738 \\ 
				& 4225 &  0.9811E-3 &1.9936& 0.1282E-2 &1.9931 \\  
				
				\hline
				$m=2$  & 36 &   0.4986E-1  & & 0.7998E-1 & \\  
				& 100 &   0.3840E-2  & 3.6986& 0.9264E-2 & 3.1099\\ 
				& 324 &0.4124E-3  & 3.2190 & 0.1103E-2 &3.0696   \\
				& 1156 &0.4904E-4 & 3.0719& 0.1362E-3 &3.0182 \\  
				
				\hline
				$m=3$	& 49 &  0.7327E-2& & 0.1476E-1 & \\  		
				& 121 &  0.6200E-3&3.5545 & 0.9512E-3 &3.9565\\  
				& 361 &0.3248E-4 & 4.2545& 0.6036E-4 & 3.9780 \\
				& 1225 &  0.1922E-5 &4.0789& 0.3765E-5 &4.0027\\  
				
				\hline
			\end{tabularx}
		\end{center}
	\end{table}

	\begin{table}[H]
		\small
		\caption{Errors and corresponding rates in the temporal direction for Example 2.} \label{tab122}
		\begin{center}
			\begin{tabularx}{1\textwidth} { 
					>{\raggedright\arraybackslash}X 
					>{\centering\arraybackslash}X 
					>{\raggedleft\arraybackslash}X     >{\raggedleft\arraybackslash}X>{\raggedleft\arraybackslash}X>{\raggedleft\arraybackslash}X  }
				\hline
				& $N$ & $E_{u}$  & Rate & $E_v$ & Rate  \\ 
				
				\hline
				$\alpha=0.4$	& 8 &  0.6592E-1  &  &0.8147E-1 & \\  
				& 16 &0.1852E-1 &1.8316 &0.2270E-1 &1.8435   \\
				& 32 & 0.4848E-2 &1.9335&0.5910E-2 & 1.9414\\  
				&64 &0.1235E-2  & 1.9724 &  0.1501E-2 & 1.9765\\	           
				\hline
				$\alpha=0.6$	& 8 &   0.6621E-1 &  &0.8173E-1 & \\  
				& 16 &0.1797E-1  & 1.8808 &0.2216E-1 &1.8825   \\
				& 32 &0.4612E-2 & 1.9628&0.5676E-2 & 1.9653\\  
				&64 & 0.1163E-2 & 1.9876 & 0.1429E-2 &1.9890\\
				\hline
				$\alpha=0.8$	& 8 &   0.6598E-1  &  &0.8152E-1 & \\  
				& 16 &0.1765E-1  & 1.9017  &0.2185E-1 & 1.8994 \\
				& 32&0.4498E-2 & 1.9726&0.5564E-2 & 1.9735 \\  
				&64 &0.1130E-2  &  1.9923 & 0.1397E-2 & 1.9929\\
				\hline
			\end{tabularx}
		\end{center}
	\end{table}

\section{Conclusions}
In this work, we solved a nonlinear time-fractional biharmonic equation using a weighted $b$-spline based mesh-free FEM with the $L2$-$1_\sigma$ approximation on a graded mesh. We derived $\alpha$-robust   \emph{a priori} bounds and convergence estimates for the proposed linearized fully-discrete scheme in the $L^2(\Omega)$ norm. Finally, we discussed two numerical experiments to validate the theoretical estimates. From the first numerical experiment, we observed that the weighted $b$-spline method provides better accuracy in comparison with the standard FEM.  \\\\

\end{document}